%Paper: funct-an/9403004
%From: Vladimir.Pestov@vuw.ac.nz
%Date: Mon, 14 Mar 1994 13:00:51 +0200

* * * * * * * * * * * * * * * * * * * * * * * * * * * * * * * * * *
* Vladimir G. Pestov         Department of Mathematics,           *
*                         Victoria University of Wellington,      *
*                            P.O. Box 600, Wellington,            *
*                                   NEW ZEALAND                   *
*                                                                 *
*  Tel. +64 4 4721-000 x 8319,   Fax +64 4 4712-070               *
*              vladimir.pestov@vuw.ac.nz                          *
* * * * * * * * * * * * * * * * * * * * * * * * * * * * * * * * * *
\documentstyle{amsppt}
\magnification=\magstep 1
\TagsOnRight
\NoBlackBoxes
\leftheadtext{V\.G\. Pestov}
% macros
\def\R {{\Bbb R}}

\def\Lie {{\operatorname{Lie}}}
\def\Card {{\operatorname {Card}}}

\def\g {{\frak g}}
\def\h {{\frak h}}
% end of macros
\topmatter
\title
Regular Lie groups and a theorem
of Lie-Palais
%$\dag$
\endtitle
\author
Vladimir G\. Pestov
\endauthor
\affil
Victoria University of Wellington, \\ P.O. Box 600, Wellington,  New
Zealand \\ \\
vladimir.pestov$\@$vuw.ac.nz
\endaffil
\abstract{In 1984 Milnor had shown how to deduce
the Lie-Palais theorem on integration of infinitesimal actions of
finite-dimensional Lie algebras
on compact manifolds from general theory of regular
Lie groups modelled on locally convex spaces.
We show how, in the case of effective action, one can
eliminate from Milnor's argument the
abstract Lie-Cartan theorem,
making the deduction rather elementary.
A machinery employed in the proof
provides a
partial solution to a problem examined in 1972 by
van Est and
\'Swierczkowski.}

\endabstract
\endsubjclass
\endkeywords
\endtopmatter
\document
%\footnote""{$\dag$ Research Report RP-94-135,
%Department of Mathematics,
%Victoria University of Wellington, January 1994.}

\subheading{1. Introduction}
A well-known result dated back to Lie
and finalized by Palais
\cite{Pa} states that every infinitesimal action of
a finite-dimensional
Lie algebra on a compact smooth manifold, $X$,
is derived from a smooth action of a
finite-dimensional Lie group on
$X$. Milnor \cite{M} gives a proof of this
theorem based on theory of regular Lie groups. This proof, however,
is only partial (in Milnor's own words), not being self-contained:
by necessity, it invokes the abstract
Lie-Cartan theorem. One definitely cannot hope to circumvent
this landmark result while proving the Lie-Palais theorem
(which simply turns into the Lie-Cartan theorem in the degenerate
case of a trivial action).
However, things are different in the most important particular
case where the infinitesimal action is
effective, that is, one deals with finite-dimensional Lie algebras
of vector fields.
We show how under this assumption the Lie-Cartan theorem
can be eliminated from the Milnor's argument, which therefore
becomes quite elementary in the sense that it does not
invoke anything advanced
beyond the scope of theory of regular
Lie groups (like the structure theory of Lie algebras).

Call a Banach-Lie algebra coming from a Lie group
{\it enlargable}. The Lie-Cartan theorem essentially says that
every finite-dimensional Lie algebra is enlargable;
in infinite dimensions non-enlargable Banach-Lie algebras
can be found \cite{vEK}. However,
a Banach-Lie algebra is enlargable as soon as it admits
a continuous monomorphism into an enlargable Banach-Lie algebra
\cite{H, \'S1, vEK}.

For more general Lie groups and algebras the situation is
much harder to deal with.
A closed Lie subalgebra of the Lie algebra of a
Fr\'echet-Lie group $G$
need not have an associated Fr\'echet-Lie group even if $G$
meets certain conditions of regularity
\cite{KYMO, M}.
For some positive advances in this direction, see \cite{L}.

We show that the above result on Banach-Lie
algebras can be pushed somewhat further to absorb a restricted version
of the Lie-Palais theorem.
Let $\h$ be the Lie algebra of a regular Lie group modelled on a
locally convex space, and let $\g$ be a Banach-Lie algebra with
finite-dimensional center
admitting a continuous monomorphism into $\h$.
Then $\g$ is enlargable.
The result is deduced from regular Lie group
theory by means of the concept of
a free Banach-Lie algebra introduced by us earlier \cite{Pe1}.
In passing, we observe that this concept
provides a partial answer to a problem investigated in 1972 by
van Est and \'Swierczkowski \cite{vE\'S}.

Since the group of diffeomorphisms of a compact manifold forms
a regular Fr\'echet-Lie group, an application of our main result
yields a theorem on integration of finite-dimensional Lie
algebras of vector fields.
The latter theorem thus merges fully in the realm of regular
Lie group theory.

\subheading{2. Preliminaries}
We first recall
a few basic facts about regular Lie groups.
A $C^\infty$ (G\^ateaux) smooth Lie group $G$ modelled on
a bornological sequentially complete locally convex space
is called {\it regular} if every smooth path
$v: I=[0,1]\to \Lie (G)$ has a left product integral
$p\colon I\to G$, that is, a solution to the differential
equation $Dp(t)=p(t)\cdot v(t)$, and if furthermore the
correspondence $v\mapsto p(0)^{-1}\cdot p(1)$ defines
a smooth map from the locally convex space
$C^\infty(I,\Lie (G))$ to $G$.
(This definition is taken from Milnor's survey \cite{M},
where regularity is understood in a somewhat less restrictive sense
than by Kobayashi {\it et al} \cite{KYMO}.
Notation is self-explanatory.)

The condition of regularity seems indispensable if one wants to
derive substantial results, and all known examples of
Lie groups modelled on locally convex spaces are regular.
Such are: the groups of diffeomorphisms of smooth compact
manifolds, certain subgroups of these,
the groups of currents (and their central extensions),
Banach-Lie groups, and many more \cite{KYMO, M}.
The following is a basic result about regular Lie groups.

\proclaim{Theorem 2.1 \text{\rm (\cite{M}; cf. \cite{KYMO})}}
Let $G$ and $H$ be smooth Lie groups
modelled on bornological sequentially complete locally convex spaces.
Let $G$ be simply connected and $H$ be regular.
Then every continuous Lie algebra morphism
$\Lie(G)\to\Lie (H)$ is tangent to a (necessarily unique)
Lie group morphism $G\to H$.
\qed\endproclaim

We also need a test for enlargability of Banach-Lie algebras.
Recall that every such algebra carries a natural structure of
a group germ \cite{\'S3}.

\proclaim{Theorem 2.2 \cite{H, \'S1}} Let $\g$ be a Banach-Lie algebra
such that there exist a group $G$ and a map
$\phi\colon\g\to G$ which is locally one-to-one and
a morphism of group germs. Then $\g$ is enlargable.
Moreover, $G$ can be given a unique
structure of a Banach-Lie group
associated to $\g$ in such a way that $\phi$ becomes an
exponential map.
\qed\endproclaim

\subheading{3. Free Banach-Lie algebras}
\proclaim{Theorem 3.1 \cite{Pe1}}
Let $E$ be a normed space. There exist a complete normed Lie algebra
$\Cal {FL}(E)$ and a linear isometry
$i_E\: E \hookrightarrow \Cal{FL}(E)$
such that:

\item{1.} The image of $i_E$ topologically generates $\Cal{FL}(E)$.

\item{2.} For every complete normed Lie algebra $\Cal L$ and an
arbitrary contracting
linear operator $f\: E \to \Cal L$ there exists a contracting
Lie algebra homomorphism $\hat f\: \Cal{FL}(E) \to \Cal L$
with $\hat f \circ i_E = f$.

The pair $(\Cal{FL}(E), i_E)$ with the  properties 1 and 2 is
unique up to an isometrical isomorphism.
\qed\endproclaim
We call $\Cal{FL}(E)$
the {\it free Banach-Lie algebra} on $E$.
It is easy to prove that if $\dim E\geq 2$ then
$\Cal{FL}(E)$ is centerless.

Let $\Gamma$ be a set. We denote the Banach-Lie algebra
$\Cal{FL}(l_1(\Gamma))$ simply by $\Cal{FL}(\Gamma)$
and call the free Banach-Lie algebra on a set $\Gamma$.
It is easy to see that the above algebra possesses the
following universal property: every map $f$ sending
$\Gamma$ to the unit ball of a Banach-Lie algebra
$\frak g$ extends in a unique way to a contracting
Lie algebra morphism from $\Cal{FL}(\Gamma)$ to $\frak g$.
The well-known fact that a Banach space of density character
$\leq\Card(\Gamma)$ is a factor space of $l_1(\Gamma)$
implies that an arbitrary Banach-Lie algebra
is a factor Banach-Lie algebra of a free Banach-Lie algebra
of the form $\Cal{FL}(\Gamma)$ \cite{Pe1}.

\definition{Remark 3.2} Not every Banach-Lie algebra
$\g$ admits a faithful representation in a
Banach space, even if $\g$ is enlargable \cite{vE\'S}.
An investigation \cite{vE\'S} of
van Est and \'Swierczkowski was partly
motivated by the problem
on whether or not every Banach-Lie algebra is a quotient algebra of a
Banach-Lie algebra with such a property.
It was suggested that
the algebra $\Lambda\,\g$ of continuous paths starting at
zero (which is enlargable and
of which $\g$ is a quotient) can be represented in
$\frak{gl}\,(E)$ for some $E$. While this conjecture remains
unproved, the following result can be considered as
one possible answer to the above motivating question.
(This fact was left unnoticed by us earlier in \cite{Pe1}.)
\enddefinition

\proclaim{Theorem 3.3} Every Banach-Lie algebra is a factor algebra of
a Banach-Lie algebra admitting a faithful representation in a
Banach space.
\endproclaim

\demo{Proof}
The left regular representation of a free Banach-Lie algebra
$\Cal{FL}(E)$ is faithful if $\dim E\geq 2$.
\qed\enddemo

\proclaim{Theorem 3.4 \cite{Pe1}}
Every free Banach-Lie algebra on a normed space is enlargable.
\endproclaim

\demo{Proof} A Banach-Lie algebra admitting a faithful
representation in a Banach space is enlargable \cite{vEK}.
\qed\enddemo

\proclaim{Corollary 3.5 \cite{\'S2}}
Every Banach-Lie algebra is a factor algebra of
an enlargable Banach-Lie algebra.
\qed\endproclaim

We shall denote by
$\Cal{FG}(\Gamma)$ the simply connected Banach-Lie group attached
to the free Banach-Lie algebra $\Cal{FL}(\Gamma)$.
It is couniversal among all
connected Banach-Lie groups of density character
$\leq\Card(\Gamma)$ \cite{Pe1}.

\subheading{4. Extension of Lie-Palais theorem in case of effective action}
\proclaim{Lemma 4.1}
Let $\g$ be a Banach-Lie algebra with finite-dimensional center
$\frak z$. Let $f$ be a continuous map from $\g$ to a
Hausdorff topological group $G$ which is a morphism of group germs.
Suppose that for every $x\in\frak z$ there is a $\lambda\in\R$
with $f(\lambda x)\neq e_G$.
Then $\g$ is enlargable.
\endproclaim

\demo{Proof} The set $\frak z\cap\ker f$ forms a closed additive
subgroup of a finite-dimensional vector space, containing no
one-dimensional linear subspaces; in other words, it is a discrete
lattice, and there exists an open neighbourhood of zero,
$U$, in $\frak z$ such that $U\cap\ker f =\{0\}$.
Choose an open $\tilde U\subset\g$ with
$0\in\tilde U\cap\frak z\subseteq U$.

Let $\phi\colon\g\to F$ be a universal morphism of a group germ
to a group. (That is, any other such morphism,
$g\colon\g\to A$, is a composition of $\phi$ and a group
homomorphism $F\to A$; cf. \cite{\'S3}.)
Let $x\in\tilde U\setminus\{0\}$;
it suffices to prove that $\phi(x)\neq e_G$
and apply Theorem 2.2.
Now, if $x\notin\frak z$, then $x$ is separated from zero by
the composition of the adjoint representation $x\mapsto
\operatorname{ad}_x$
and the exponentiation $\operatorname{End}(\frak g_+)\to
\operatorname{GL}(\g_+)$ (both being morphisms of group germs).
If $x\in\frak z$, then $x\in U$ and $f(x)\neq e_G$.
\qed\enddemo

\proclaim{Main Theorem 4.2} Let $\g$ be a Banach-Lie algebra
with finite-dimensional center
and let $\g$ admit a continuous monomorphism into the Lie algebra
of a regular Lie group. Then $\g$ is enlargable.
\endproclaim

\demo{Proof} Let $H$ be a regular Lie group
 and
$h\colon\g\to\Lie(H)$ be a continuous Lie monomorphism.
Denote by $B$ the unit ball in $\g$, and by $i$ the
contracting Lie algebra morphism
$\Cal{FL}(B)\to\g$ extending the identity map
$B\to B$.  According to Theorem 2.1,
the composition morphism $h\circ i\colon \Cal{FL}(B)\to\Lie(H)$
is tangent to a morphism of Lie groups
$\hat i\colon\Cal{FG}(B)\to H$.
Denote $I=\ker h\circ i$, $J=\ker\hat i$, and let $A$ be
a subgroup of $\Cal{FG}(B)$ algebraically generated by the
exponential image of $I$.
Let $G$ denote the Hausdorff topological
group $\Cal{FG}(B)/J$, and let
$f\colon \g\to G$ be a map defined by
$f(x+I)=(\exp x)\cdot J$. Since $A\subseteq J$ and
$\exp_{\Cal{FG}(B)}$ is a morphism of group germs,
the map $f$ is correctly defined and is
a morphism of group germs;
diagrammatic search shows that it is continuous.
The canonical group monomorphism $\psi\colon G\to H$ is
continuous as well, and
$\psi\circ\phi =\exp_H\circ h$.
Since for every $x\in\Lie (H)$ one surely has
$\exp_H(\lambda x)\neq e_H$ for an appropriate $\lambda\in\R$,
the triple $(\g,G,f)$ is under conditions of
Lemma 4.2.
\qed\enddemo

\proclaim{Corollary 4.3}
Every continuous
monomorphism from a Banach-Lie algebra with finite-dimensional
center to the Lie algebra
of a regular Lie group is tangent to an appropriate Lie group morphism.
\endproclaim

\demo{Proof}
Results from direct application of Theorems 4.2 and 2.1.
\qed\enddemo

As a corollary, we can now prove a restricted version
of the Lie-Palais theorem by methods of regular Lie group theory.

\proclaim{Theorem 4.4 \text{\rm (Lie-Palais)}}
Let $M$ be a compact closed
manifold and let $\g$ be a finite-dimensional
Lie algebra. Then every effective
infinitesimal action of $\g$ on $M$ is derived from
a smooth global action of a finite-dimensional Lie group on $M$.
\endproclaim

\demo{Proof} The effective
infinitesimal action of $\g$ upon $M$ can be
viewed as a Lie algebra monomorphism from $\g$ to the algebra
$\operatorname{vect}M$ of smooth vector fields on $M$ endowed with
the $C^\infty$-topology. Clearly, the conditions of the Main Theorem
4.2 and Corollary 4.3
are satisfied, and the resulting Lie group morphism from
a Lie group $G$ attached to $\g$ to
the diffeomorphism group $\operatorname{diff}M$
(which is a regular Lie group assigned to $\operatorname{vect}M$,
\cite{KYMO, M}) determines a desired global smooth action.
\qed\enddemo

\definition{Concluding Remarks 4.5. 1}
We do not know whether the condition upon the
center being finite-dimensional can be relaxed in the
Main Theorem 4.2.
However, our Example in \cite{Pe2} can be remade easily to show that
the Lemma 4.1 is no longer valid for Banach-Lie algebras with
infinite-dimensional center.
\smallskip
\noindent{\bf 2.} Does not the present result shed a
ray of hope, wee as it is,
on an (apparently, abandoned) search for an ``easy''
proof of the Lie-Cartan theorem itself?
(See \cite{S}, II.V.8, and \cite{Po, Pe2}.)
\smallskip
\noindent{\bf 3.} Another reputed open problem,
which might well fall within this circle of ideas,
is that of
existence of an effective transitive action of
a (connected) infinite-dimensional
Banach-Lie group on a compact manifold
\cite{O, OdlH}.
\enddefinition

\subheading{Acknowledgments}
It is appropriate to thank Joshua Leslie for a
few encouraging and stimulating
letters sent by him in 1987 to the author who then
worked at a secluded Siberian university.
A Victoria University 1992-93 Small Research Grant is acknowledged.

\bye